\documentclass[12pt]{amsart}

\usepackage{amssymb} 

\theoremstyle{plain}
\newtheorem{theorem}{Theorem}[section]
\newtheorem{proposition}{Proposition}[section]
\newtheorem{corollary}{Corollary}[section]
\newtheorem{lemma}{Lemma}[section]

\theoremstyle{definition}

\newtheorem{example}{Example}[section]
\newtheorem{remark}{Remark}[section]

\newcommand{\sn}{{\mathbb S}}
\newcommand{\rn}{{\mathbb R}}

\DeclareMathOperator{\trace}{trace}

\DeclareMathOperator{\ricci}{Ricci}
\DeclareMathOperator{\Div}{div}

\begin{document}
\title{Harmonic maps and sections on spheres}

\author{M. Benyounes}
\address{D{\'e}partement de Math{\'e}matiques, UMR 6205 \\
Universit{\'e} de Bretagne Occidentale \\
6, avenue Victor Le Gorgeu \\
CS 93837, 29238 Brest Cedex 3, France}
\email{Michele.Benyounes@univ-brest.fr {\rm{and}} Eric.Loubeau@univ-brest.fr}

\author{E. Loubeau}

\author{C.~M. Wood}
\address{Department of Mathematics \\
University of York \\
Heslington, York Y010 5DD, U.K.}
\email{cmw4@york.ac.uk}


\date{\today}

\begin{abstract}
The absence of interesting harmonic sections for the Sasaki and Cheeger-Gromoll metrics has led to the consideration of alternatives, 
for example in the form of a two-parameter family of natural metrics shown to relax existence conditions for harmonicity~\cite{LBW}.\\
This article investigates harmonic Killing vector fields, proves their non-existence on ${\sn}^{2}$, 
obtains rigidity results for harmonic gradient vector fields on the two-sphere, classifies spherical quadratic gradient fields 
in all dimensions and determines the tension field, concluding with the discovery of a family of metrics making 
Hopf vector fields harmonic maps on ${\sn}^{2n+1}$.
\end{abstract}

\maketitle

\section{Introduction}

Since the tangent bundle $TM$ of a manifold $M^n$ can be equipped with the structure of a 
$2n$-dimensional manifold, it admits, at each of its points $e$, a tangent space of 
its own, denoted $T_{e}TM$, their disjoint union making up the bitangent bundle $TTM$.
\newline For a Riemannian manifold $(M,g)$, the differential of the canonical projection $\pi :TM\to M$ and the Levi-Civita connection 
$\nabla$, via the connection 
map $K: TTM \to TM$, induce a natural splitting of $TTM$ 
into vertical and horizontal subbundles so that:
$$ T_{e}TM = V_{e} \oplus H_{e} .$$
As the restriction of $K$ to $V_{e}$ leads to a canonical identification with $T_{x}M$ ($\pi(e)=x$), we can lift vector fields on $M$ to 
vertical 
and horizontal vector fields on $V$ and $H$ (cf.~\cite{Gud-Kap-Expo} for details).
\newline The simplest metric one can define on $TM$ is the Sasaki metric, which requires the projections from $V$ and $H$ to $TM$ to 
be isometries and the two bundles to be orthogonal. Unfortunately, at least on compact manifolds, this most natural of metrics does not 
allow the 
existence of interesting harmonic maps or even sections (see definitions below). In order to obtain some flexibility, a new 
(two-parameter) family of metrics on $TM$ was introduced in~\cite{LBW}:
$$h_{p,q}(A,B) = g(d\pi(A),d\pi(B)) + \omega^{p}(e)\left[ g(KA,KB) + qg(KA,e)g(KB,e)\right] ,$$
where, $e\in TM$, $A,B \in T_{e}TM$, $\omega(e) = (1+|e|^2)^{-1}$ and $p,q \in \rn$.
\newline While it is possible to allow $q$ negative, $h_{p,q}$ being then only Riemannian in a tubular neighbourhood of the zero 
section, we will assume $q\geqslant 0$, except for a few special cases.
\newline Note that $h_{0,0}$ is precisely the Sasaki metric and $h_{1,1}$ the Cheeger-Gromoll metric (cf.~\cite{Kowal,Gud-Kap-Tok} for 
more on 
their geometries).
\newline These metrics $h_{p,q}$, named ``generalised Cheeger-Gromoll metrics'', are all natural metrics (i.e. $d\pi$ is Riemannian 
and $V$ and $H$ are orthogonal). Their geometry is studied in~\cite{geom}.

\section{$(p,q)$-harmonic sections}
Assuming the base compact (otherwise work with compact subsets), we consider, rather than the ``full'' Dirichlet energy 
functional, the vertical energy of a vector field $\sigma$:
$$ E^{v}(\sigma) = \tfrac{1}{2} \int_{M} |d^{v}\sigma|^2   v_{g}$$
(these two only differ by an additive constant) and derive its Euler-Lagrange equation for 
variations through sections~\cite{LBW}:
\begin{equation}\label{eq1}
(1+|\sigma|^2) \nabla^{*}\nabla \sigma + 2p \nabla_{X(\sigma)}\sigma = \left[ p|\nabla\sigma|^{2} -pq |X(\sigma)|^2 
-q(1+|\sigma|^2)\Delta(|\sigma|^2/2) \right] \sigma
\end{equation}
where $\nabla^{*}\nabla \sigma = - \trace(\nabla^{2})\sigma$ (same sign convention for the Laplacian on functions) and $X(\sigma) = 
\nabla\left(\tfrac{|\sigma|^2}{2}\right)= g(\nabla\sigma,\sigma)$.
\newline Clearly, parallel vector fields (when they exist) are $(p,q)$-harmonic sections and the first non-parallel examples are the 
standard Hopf vector fields on ${\sn}^{2p+1}$ ($p=2$ and any $q$). Besides, rescaling them yields examples for all $(p,q)$ with 
$p>1$. Moreover, violating our convention $q\geqslant 0$, conformal vector fields on $\sn^n$ are $(p,q)$-harmonic if and only if 
$n\geqslant 3$, 
$p=n+1$ and $q=2-n$ (\cite{LBW}).
\newline On the other hand, when the base is compact (without boundary) and $p\leqslant 1$ ($q\geqslant 0$), then a vector field is 
$(p,q)$-harmonic if and only if it is parallel (cf.~\cite{LBW}).
\newline This includes the case of the original Cheeger-Gromoll metric ($p=q=1$) (see also~\cite{Cez}).
Note our convention for the curvature tensor: $ R(X,Y) = [\nabla_{X}, \nabla_{Y}] - \nabla_{[X,Y]} .$

\section{Killing vector fields}

As illustrated by Hopf fields on $\sn^{2p+1}$, Killing vector fields are prime candidates for $(p,q)$-harmonicity.

\begin{proposition}
A Killing vector field $\sigma$ is $(p,q)$-harmonic if and only if:
\begin{equation}\label{eq2}
(1+ |\sigma|^{2}) \ricci(\sigma) - 2p \nabla_{\nabla_{\sigma}\sigma}\sigma = \Big[ 
p|\nabla \sigma|^{2} - pq |\nabla_{\sigma}\sigma|^{2} - 
q(1+ |\sigma|^{2})\Delta(\tfrac{1}{2}|\sigma|^{2}) \Big] \sigma
\end{equation}
\end{proposition}
\begin{proof}
For $\sigma$ Killing, $X(\sigma) =  - \nabla_{\sigma} \sigma$ and $\nabla^{2}_{X,Y} \sigma = - R(\sigma,X)Y$ so 
$\nabla^{*}\nabla \sigma = \ricci(\sigma)$.
\end{proof}

The first sub-case to investigate is constant norm.

\begin{proposition} \label{prop2.2}
A Killing vector field $\sigma$ of constant norm $|\sigma|^2 = k^2$ is $(p,q)$-harmonic if and only if 
it is an eigenvector of the Ricci operator. Then $p= 1 + \tfrac{1}{k^2}$ and $q$ is any (positive) real number.
\end{proposition}
\begin{proof}
If $\sigma$ is Killing and $|\sigma|^2 = k^2$, then $X(\sigma) =0$ and $\ricci(\sigma,\sigma)=|\nabla \sigma|^2$ (\cite{Petersen}), so 
$(p,q)$-harmonicity becomes:
$$(1+ k^{2}) \ricci(\sigma) = ( p|\nabla \sigma|^{2}) \sigma 
=  p\ricci(\sigma,\sigma) \sigma.$$ 
\end{proof}

\begin{example}\label{example1}
There are easy examples of Killing vector fields which are eigenvectors of 
Ricci on the space forms ${\rn}^n$ and ${\sn}^{2p+1}$. In dimension three, one can 
also search among Thurston's geometries, starting with ${\sn}^{2}\times {\rn}$ and 
$\tfrac{\partial}{\partial t}$. Unfortunately, all these examples, except Hopf fields, are parallel.
\newline To obtain new non-parallel examples, we turn to the Heisenberg space $H_{3}$, 
identified with $\rn^3$ endowed with the metric $dx^2 + dy^2 + (dz -xdy)^2$ and 
$\widetilde{SL_{2}(\rn)}$, the universal cover of the special linear group, seen as 
$\rn^3_{+}=\{(x,y,z)\in\rn^3 : z > 0\}$ with the metric 
$(dx+\tfrac{dy}{z})^2+\tfrac{dy^2+dz^2}{z^2}$. 
\newline For $H_{3}$, the vectors $E_{1}= \tfrac{\partial}{\partial x},  
E_{2}=\tfrac{\partial}{\partial y}+ x \tfrac{\partial}{\partial z}, 
E_{3}=\tfrac{\partial}{\partial z}$ form an orthonormal basis and $E_{3}$ is Killing 
and non-parallel. Besides one verifies that 
$\ricci(E_{3}) = \tfrac{1}{2} E_{3}$.
\newline For $\widetilde{SL_{2}(\rn)}$, things are similar, since an orthonormal basis 
is given by $E_{1}= \tfrac{\partial}{\partial x}$, 
$E_{2}=z\tfrac{\partial}{\partial y}- \tfrac{\partial}{\partial x}$, 
$E_{3}=z\tfrac{\partial}{\partial z}$ and $E_{1}$ is Killing and non-parallel, moreover $\ricci(E_{1}) = \tfrac{1}{2} E_{1}$.
\end{example}

More generally, constancy of the norm is replaced with harmonicity.

\begin{proposition}\label{prop1}
Let $\sigma$ be a $(p,q)$-harmonic Killing vector field with harmonic norm.
If $p\geqslant 0$ and $(p-1)|\sigma|^2 \leqslant 1$ (automatic when $p\leqslant 1$) then $\sigma$ is of constant norm.

\end{proposition}

\begin{proof}
In this case, $\ricci(\sigma,\sigma) = 
|\nabla\sigma|^{2}$ and $\sigma$ $(p,q)$-harmonic implies:
$$ \left( 1 + (1-p)|\sigma|^{2} \right) \ricci(\sigma,\sigma) = 
-p (2 + q|\sigma|^{2})|\nabla_{\sigma} \sigma|^2 $$
hence $\sigma$ has constant norm.
\end{proof}

\begin{remark}
On compact manifolds, Proposition~\ref{prop1} is only a minute improvement on~\cite[Theorem~4.6]{LBW}.
\end{remark}

\begin{remark}
Non-parallel Killing vector fields do not exist on a compact manifold with negative $\ricci$. On the 
other hand, when $\ricci > 0$, if $|\sigma|^2$ is harmonic then it must be constant~\cite{Yau}.
\end{remark}

The two-sphere is a particularly restrictive case.

\begin{proposition}\label{prop3.4}
There is no non-zero $(p,q)$-harmonic Killing vector field on $\sn^2$ ($p > 0$, $q\geqslant 0$).
\end{proposition}
We first need:
\begin{lemma}
Let $\sigma$ be a non-zero Killing vector field on $\sn^2$. Then:\\
i) The set $A = \{ x \in {\sn^2} : \sigma(x) =0 \}$ has empty interior.\\
ii) If $\nabla_{\sigma} \sigma = 0$ on an open subset $U$ then 
$\Delta \tfrac{|\sigma|^2}{2} = 0$ on $U$.
\end{lemma}
\begin{proof}
i) If $\sigma\equiv 0$ on an open subset $B$ of $\sn^2$, then 
$(\nabla_{X} \sigma)(x) = 0, \forall x\in B , \forall X \in T_{x}{\sn^2}$ and $\sigma \equiv 0$ on $\sn^2$ (\cite{Petersen}).\\
ii) Assume $\nabla_{\sigma} \sigma \equiv 0$ on an open subset $U$ of $\sn^2$, then $|\sigma|^2$ is constant and $\Delta 
\tfrac{|\sigma|^2}{2} = 0$ on $U$.
\end{proof}
\begin{proof}[Proof of Proposition~\ref{prop3.4}]
Let $\sigma$ be a non-zero $(p,q)$-harmonic Killing vector field on $\sn^2$, $U$ and $A$ as in the above lemma.\\
On $W = {\sn^2} \setminus (U \cup A)$, $\sigma /|\sigma|$ and $\nabla_{\sigma} 
\sigma / |\nabla_{\sigma} \sigma|$ form an orthonormal basis, put 
$V= \nabla_{\sigma} \sigma$ then:\\
$$ |\nabla \sigma|^2 = \tfrac{1}{|\sigma|^2} |\nabla_{\sigma} \sigma|^2 + 
\tfrac{1}{|V|^2} |\nabla_{V} \sigma|^2 .$$
Since $\sigma$ is $(p,q)$-harmonic:
$$\nabla_{\nabla_{\sigma}\sigma}\sigma = \tfrac{-1}{2p} \Big( p|\nabla \sigma|^{2} - 
pq |\nabla_{\sigma}\sigma|^{2} - q(1+|\sigma|^{2})\Delta(\tfrac{1}{2}|\sigma|^{2}) 
- (1+ |\sigma|^{2}) \Big) \sigma $$
so that $\nabla_{\nabla_{\sigma}\sigma}\sigma = f 
\sigma$ with $f = \tfrac{-1}{2p} \Big( p|\nabla \sigma|^{2} - pq |\nabla_{\sigma}\sigma|^{2} - 
q(1+|\sigma|^{2})\Delta(\tfrac{1}{2}|\sigma|^{2}) 
- (1+ |\sigma|^{2}) \Big)$.
On $W$:
$$|\nabla_{\sigma} \sigma|^2 = -f |\sigma|^2 \quad \mbox{and } |\nabla_{V} \sigma|^2 = f^2 |\sigma|^2 .$$
Therefore $ |\nabla \sigma|^2 = -2 f.$
Note that by continuity this holds on $\sn^2 \setminus U$ and $f$ cannot vanish.\\
Hence
$$ 0 =  pq |\nabla_{\sigma}\sigma|^{2} + q(1+|\sigma|^{2}) 
\Delta(\tfrac{1}{2}|\sigma|^{2}) + (1+ |\sigma|^{2})$$ 
so that $|\sigma|^2$ is constant.
\end{proof}

However, one can try to obtain $(p,q)$-harmonic vector fields by multiplying a Killing vector field by a function.

\begin{example}
Consider $\sigma(x,y,z)=(-y,x,0)$ the Killing vector field on $\sn^2$ given by rotation around the $z$-axis and complete it into an 
orthonormal frame $\{\tfrac{\sigma}{x^2+y^2} , \tfrac{\gamma}{x^2+y^2} \}$ of $\sn^2 \setminus \{ N,S\}$, the sphere less the poles.\\
For $f\in C^{\infty}(\sn^2)$, let $V= f\sigma$, then standard computations show that $V$ $(p,q)$-harmonic requires 
$\sigma(f)=0$, i.e. $f(x,y,z) = F(x^2 + y^2)$ for some real function $F$.
\newline Then $V= F(x^2 + y^2)\sigma$ is $(p,0)$-harmonic if:
\begin{equation}\label{eq*}
(1 + tF^{2})\big[ 4t(t - 1)F'' + (10t-8) F' +F \big] + 4p(1-t)\big[ tF^2 F' + t^{2}F(F')^{2}\big] =0
\end{equation}
Clearly constant functions are not solutions and while $F(t) = \tfrac{1}{\sqrt{(p-1)t}}$ is a solution, it does not extend over the 
poles.
\newline More generally, if $F$ is a solution of \eqref{eq*} on $[0,1]$, assuming (wlog) $F(0)\geqslant 0$ and $F'(t) \geqslant 0$ for 
$t\in 
]0,\epsilon[$ ($\epsilon > 0$), then $F' \geqslant \tfrac{1}{\sqrt{t^{4}(1-t)}}$ on $]0,\epsilon[$. 
Hence no function $F$ on $\sn^2$ can make $F\sigma$ $(p,0)$-harmonic.
\end{example}

\section{Gradient vector fields}

First, we find an obstruction for another important class of vector fields on $\sn^2$.

\begin{proposition}\label{proposition2}
Let $\sigma$ be a (non-zero) $(p,q)$-harmonic vector field defined by $\sigma = \nabla f$ with $\Delta 
f = \lambda f$. If $\nabla_{\sigma}\sigma$ and $\sigma$ are non-colinear on an open 
subset $U$ of $\sn^2$, then $\Delta\tfrac{|\sigma|^2}{2} \leqslant 0$ on $U$.
\end{proposition}

\begin{proof}
From standard properties of gradient vector fields, $\sigma = \nabla f$ with $\Delta 
f = \lambda f$, is $(p,q)$-harmonic if:
$$(1+ |\sigma|^2)(\lambda -1)\sigma + 2p \nabla_{\nabla_{\sigma}\sigma}\sigma = \left[ 
p|\nabla\sigma|^2 - pq |\nabla_{\sigma}\sigma|^2 - 
q(1+|\sigma|^2) \Delta \tfrac{|\sigma|^2}{2}\right] \sigma$$
Therefore $\nabla_{\nabla_{\sigma}\sigma}\sigma = f \sigma$ with: 
$$f = \tfrac{1}{2p} \left( p|\nabla\sigma|^2 - pq |\nabla_{\sigma}\sigma|^2 - 
q(1+|\sigma|^2) \Delta \tfrac{|\sigma|^2}{2} - (\lambda -1)(1+ |\sigma|^2)\right).$$
Note that $\sigma$ cannot vanish on an open subset of $\sn^2$. Assume that $\left\{ \tfrac{\sigma}{|\sigma|} , \tfrac{b}{|b|} 
\right\}$, where $b= 
\nabla_{\sigma} \sigma - \langle \nabla_{\sigma} \sigma , \sigma\rangle  \tfrac{\sigma}{|\sigma|^2}$, 
is an orthonormal basis on an open subset $U$, then:
$$|\nabla \sigma|^2 = |\nabla_{\tfrac{\sigma}{|\sigma|}} \sigma|^2 + 
|\nabla_{\tfrac{b}{|b|}} \sigma|^2 $$
with $$|\nabla_{\tfrac{\sigma}{|\sigma|}} \sigma|^2 =  
\tfrac{1}{|\sigma|^2} \langle f\sigma , \sigma \rangle  = f$$
while $$|\nabla_{\tfrac{b}{|b|}} \sigma|^2 = \tfrac{1}{|b|^2} |\nabla_{b} \sigma|^2 $$
$$|b|^2 = f|\sigma|^2 - \tfrac{1}{|\sigma|^2}\langle \nabla_{\sigma} \sigma, \sigma\rangle ^2
\quad \mbox{and }
|\nabla_{b} \sigma|^2 = f\left[ f|\sigma|^2 - \tfrac{1}{|\sigma|^2}\langle \nabla_{\sigma} \sigma , \sigma\rangle ^2  
\right] $$
Therefore $\tfrac{1}{|b|^2} |\nabla_{b} \sigma|^2 = f$ and $|\nabla \sigma|^2 = 2f$.\\
Plugging in the definition of $f$, we have:
$$ 0= -pq|\nabla_{\sigma}\sigma|^2 -q(1+|\sigma|^2) \Delta\tfrac{|\sigma|^2}{2} - 
(\lambda -1)(1+|\sigma|^2),$$
therefore $\Delta\tfrac{|\sigma|^2}{2} \leqslant 0$ on $U$.\\
If $U$ is dense in $\sn^2$, this forces $|\sigma|^2$ to be constant.
\end{proof}

\begin{remark} For conformal vector fields on $\sn^{n}$, $\nabla_{\sigma}\sigma$ and $\sigma$ are colinear. 
However \cite[Theorem 5.4]{LBW} states that conformal vector fields on $\sn^2$ cannot be $(p,q)$-harmonic.
\end{remark}

\subsection{Spherical Quadratic Gradient Fields}

To a symmetric bilinear form $\beta$ on ${\rn}^{n+1}$, associate the linear map $B$:
$$
B(x)\centerdot y = \beta(x,y),
\quad\text{for all $x,y\in{\rn}^{n+1}$,}
$$
define the restriction of the associated quadratic form to
${\sn}^n$:
$$
\lambda(x) =\beta(x,x) = Bx\centerdot x,
$$
and $\sigma=\tfrac{1}{2}\nabla\lambda$ its spherical gradient. Clearly $\sigma(x) = B(x)-\lambda(x)x$.
More generally, for any $k\in\mathbb N$: 
$$\lambda_k(x) = B^k(x)\centerdot x, \quad \mbox{and } \sigma_k(x) = \tfrac12 \nabla\lambda_k = B^k(x)-\lambda_k(x) x.$$
Consider the symmetric $(1,1)$-tensors on ${\sn}^n$:
$$L_kX = B^kX-(B^kX\centerdot x) x,$$ 
and put $L=L_1$. Note that: $L_k\neq L^k$ if $k>1$.

By direct calculations, we obtain:
\begin{proposition}\label{cprop1}
\begin{align*}
&{(a)}\quad \langle \sigma_k,\sigma_l\rangle  = \lambda_{k+l}-\lambda_k \lambda_l \\
&{(b)}\quad L_l \sigma_k = \sigma_{k+l}-\lambda_k \sigma_l \\
&{(c)}\quad \nabla_{X}\sigma_k = L_k X - \lambda_k X \\
&{(d)}\quad \Div\sigma_k = \trace(B^k) - (n+1)\lambda_k \\
&{(e)}\quad | L_k |^2 = | B^k |^2 - 2\lambda_{2k} + \lambda_k^{ 2} \\
&{(f)}\quad \nabla_{X}L_k(Y) = -\langle X,Y \rangle \sigma_k - \langle \sigma_k,Y\rangle X \\
&{(g)}\quad \nabla^{2}_{X,Y}\sigma_k = -\langle X,Y\rangle \sigma_k - \langle \sigma_k,Y\rangle X - 2\langle \sigma_k,X\rangle Y
\end{align*}
\end{proposition}

\begin{corollary}\label{ccol2}
 {(a)}\quad
When $L_k$ is viewed as a $TM$-valued $1$-form on $M$: $\quad\delta L_k = (n+1) \sigma_k$
\newline {(b)}\quad
$\nabla^*\nabla\sigma_k = (n+3)\sigma_k$
\end{corollary}

A straightforward consequence of Proposition~\ref{cprop1} is:
\begin{proposition}\label{cprop3}
For $\sigma=\tfrac{1}{2}\nabla\lambda$, we have:
\begin{align*}
&{(a)}\quad |\sigma|^2 = \lambda_2-\lambda^2 \\
&{(b)}\quad X(\sigma) = \sigma_2-2\lambda \sigma \\
&{(c)}\quad \nabla_{X}\sigma = LX-\lambda X \\
&{(d)}\quad \nabla_{X(\sigma)}\sigma = \sigma_3 - 3\lambda \sigma_2 + (4\lambda^2-\lambda_2)\sigma \\
&{(e)}\quad | \nabla\sigma |^2 = | B |^2 - 2\lambda_2 - 2\lambda\trace B 
+ (n+3)\lambda^2 \\
&{(f)}\quad \nabla^*\nabla\sigma = (n+3)\sigma \\
&{(g)}\quad \Delta |\sigma|^2 = -2 | B |^2 + 2(n+5) \lambda_2 + 4\lambda\trace B 
- 4(n+3)\lambda^2
\end{align*}
\end{proposition}

This formalism leads to a classification.

\begin{theorem}\label{thclass}
A quadratic gradient vector field $\sigma$ on $\sn^n$ is $(p,q)$-harmonic if and only if $n$ is odd, $n\geqslant 5$, $p=\tfrac{n+3}{2}$ 
and $q$ 
is the only solution (necessarily negative) of:
$$ 8(n+1) q^2 +2(3n^2 -2n -17) q + (n-3)(n^2 -9) =0$$
such that $\mu^2 = \tfrac{9-n^2}{4q} -2(n+3)$ is strictly positive.
\newline Besides, $\sigma$ is given by $\tfrac{1}{2}\nabla \lambda$ where $\lambda(x) = Bx.x$ and $B$ is a symmetric matrix 
with two eigenvalues $\alpha$ and $\beta$, of equal multiplicities $\tfrac{n+1}{2}$, such that $|\alpha - \beta| = \mu$.
\end{theorem}

\begin{corollary}\label{rmkA}
There is no $(p,q)$-harmonic quadratic gradient vector fields on even-dimensional spheres, and, more generally, none for 
positive values of $q$ (recall that $h_{p,q}$ is a Riemannian metric on the whole of $T\sn^n$ if and only if $q$ is positive).
\end{corollary}

\begin{remark}\label{rmkB}
a) $B$ and $\tilde{B} = B - \gamma Id$ give rise to the same gradient vector field.
\newline b) Since $B$ can be written $B= 
P^{-1} A P$, $A$ diagonal and $P$ orthogonal, $\lambda(x) = B(x).x$ and 
$\tilde{\lambda}(x) = A(x).x$ are isometric and so are $\sigma = \tfrac{1}{2} \nabla \lambda$ 
and $\tilde{\sigma} = \tfrac{1}{2} \nabla \tilde{\lambda} $. More prosaically:
$$\lambda(x) = \langle B(x),x\rangle  = \langle P^{-1} A P x,x\rangle = \langle Ax,x\rangle \circ\phi ,$$
for some isometry $\phi$.
\newline So, up to isometries and constants, the $(p,q)$-harmonic vector fields 
of Theorem~\ref{thclass} are the gradients of $\lambda(x) = x_{1}^2 + \cdots + x_{\tfrac{n+1}{2}}^2$.
\end{remark}

\begin{proof}[Proof of Theorem~\ref{thclass}]

From Proposition~\ref{cprop3}:
\begin{align} \label{eq1a}
&(1+|\sigma|^2) \nabla^*\nabla\sigma + 2p \nabla_{X(\sigma)}\sigma =\\
& 2p\sigma_3 - 6p\lambda\sigma_2 
+ \bigl(n+3 + (n+3-2p)\lambda_2 + (8p-n-3)\lambda^2\bigr)\sigma \notag
\end{align}
hence, if $\sigma$ is $(p,q)$-harmonic, $\sigma_3 - 3\lambda\sigma_2$ must be colinear to $\sigma$.

\begin{remark}
a) It is easy to see that $\sigma_2=m\sigma$ for $m\in{\rn}$ if and
only if $B$ satisfies $B^2-mB=cI$ for $c\in{\rn}$.  \\
b) $B^2-mB=cI$ is equivalent to $B$ having either one or two
distinct eigenvalues, depending on whether the roots of the polynomial
$x^2-mx-c$ are repeated or distinct. If $B$ has a single eigenvalue
$\mu$ then $m=0$ and $c=\mu^2$, whereas if $B$ has distinct eigenvalues
$\mu,\nu$ then $m=\mu+\nu$ and $c=-\mu\nu$, and:
\begin{align*}
\sigma(x) &= (\mu-\nu)| x_\nu |^2 x_\mu 
+ (\nu-\mu)| x_\mu |^2 x_\nu \\
\sigma_2(x) &= (\mu^2-\nu^2)| x_\nu |^2 x_\mu 
+ (\nu^2-\mu^2)| x_\mu |^2 x_\nu,
\end{align*}
where $x=x_\mu+x_\nu$ is the decomposition into eigenvectors.
\end{remark}

\begin{lemma}
A non-trivial $\sigma$ satisfies:
\begin{equation}
\sigma_3 - 3\lambda\sigma_2 = f\sigma,
\label{cmmidstar}
\end{equation} 
for some smooth function $f\colon {\sn}^n\to{\rn}$ if and only if $B$ has
precisely two distinct eigenvalues.
\end{lemma}

\begin{proof}
If $B$ has two distinct eigenvalues $\mu,\nu$, then $B^2-mB=cI$, where $m=\mu+\nu$ and $c=-\mu\nu$ and:
$$
\sigma_3-3\lambda\sigma_2 = (m^2-3m\lambda+c)\sigma
$$
To show the converse, let $\{e_i\}$ be an orthonormal $B$-eigenbasis of
${\rn}^{n+1}$:
$$
B e_i = \mu_i e_i,\quad
\text{where $\mu_i\in{\rn}$.}
$$
Notice that for all $k\in\mathbb N$ we have:
\begin{align*}
&\lambda_k(e_i) =  \mu_i^{ k} \quad \mbox{and } \sigma_k(e_i) = 0 ,
\end{align*}
then $x=e_i\in {\sn}^n$ and covariantly differentiating Equation~\eqref{cmmidstar}
along $X=e_j \in T_x{\sn}^n$ for some $j\neq i$ yields:
$$
\nabla_{X}\sigma_3 - 3\mu_i \nabla_{X}\sigma_2 = f(e_i) \nabla_{X}\sigma
$$
which reduces to:
\begin{equation*}
\mu_j^{ 3}-\mu_i^{ 3} - 3\mu_i(\mu_j^{ 2}-\mu_i^{ 2}) 
= (\mu_j-\mu_i) f(e_i),
\quad\text{for all $i\neq j$.}
\end{equation*}
If $\sigma$ is non-trivial, there exist $i,j$
such that $\mu_i\neq\mu_j$, and:
\begin{align}
f(e_i) &= \mu_j^{ 2} - 2 \mu_i \mu_j - 2 \mu_i^{ 2} 
\label{cddag} \\
f(e_j) &= \mu_i^{ 2} - 2 \mu_i \mu_j - 2 \mu_j^{ 2}
\label{cdagdag}
\end{align}
If $\mu_k$ is a third distinct eigenvalue then:
$$
f(e_i) = \mu_k^{ 2} - 2 \mu_i\mu_k - 2\mu_i^{ 2}
$$
and comparing with \eqref{cddag} and \eqref{cdagdag} yields:
$$\mu_i = \frac12(\mu_j+\mu_k) \quad \mbox{and }
\mu_j = \frac12(\mu_i+\mu_k) .$$
Therefore:
$$
\mu_j = \frac12(\mu_i+\mu_k) = \frac14(\mu_j+3\mu_k)
$$
hence $\mu_j=\mu_k$. So $B$ has at most two distinct eigenvalues.
\end{proof}

\noindent{\em End of proof of Theorem~\ref{thclass}:}\\
By Remark~\ref{rmkB} a), one eigenvalue, $\nu$, can be chosen zero, and the other, $\mu$, strictly positive, of multiplicity $k$. 
Decomposing $x = x_{\mu} + x_{\nu} \in \sn^n$ into eigenvectors, the left-hand side of \eqref{eq1a} becomes:
\begin{align}
2p\mu^2 + n+3 +( n+3-8p) \mu^2 |x_{\mu}|^2 + (8p -n-3) \mu^2 |x_{\mu}|^4 \label{eq2a}
\end{align}
while its right-hand side gives:
\begin{align}
&k\mu^2(p+q) + \mu^2[-2p(k+1) -pq\mu^2 +kq\mu^2 -q(n+5+2k)] |x_{\mu}|^2 + \label{eq3a}\\
& \mu^2[p(n+3) + 5pq\mu^2 -q(n+5+2k)\mu^2 -kq\mu^2 + 2(n+3)q]|x_{\mu}|^4 + \notag\\
&q\mu^4[-8p + 2(n+3) + n+5+2k]|x_{\mu}|^6 + 2q\mu^4(2p- n-3)|x_{\mu}|^8 \notag
\end{align}

Note that by continuity, \eqref{eq2a} and \eqref{eq3a} must be equal everywhere. Equating the coefficients yields $p=\tfrac{n+3}{2}$,  
$k=\tfrac{n+1}{2}$ (but rules out $q=0$) and
the system:
\begin{equation}
\left\lbrace
\begin{array}{l}
\tfrac{n^2 -9}{4} + q(\mu^2 + 2(n+3)) =0  \quad \mbox{ (forcing $q<0$ and excluding $n=3$)}\\
8(n+1) q^2 +2(3n^2 -2n -17) q + (n-3)(n^2 -9) =0
\end{array}
\right.
\end{equation}
of the form:
\begin{equation}
\left\lbrace
\begin{array}{l}
(\mu^2 + a)q = K  \\
\mu^2 (q +b) = K'
\end{array}
\right.
\end{equation}
with $K$ negative and $K'$ positive. These represent two (branches) of hyperbolas in different quadrants  but intersecting (once and 
only once).
\end{proof}

\section{Harmonic maps}
With the same energy functional $E^v$ on sections, one can formulate several variational problems, according to the 
type of variations allowed. The lifting of all restrictions, that is permitting variations through mere maps, 
brings us back to harmonic maps (cf.~\cite{E-L}). In this situation, the associated Euler-Lagrange operator is the tension field 
$\tau(\sigma) = 
\trace \nabla d\sigma$.\\
Standard computations provide:
\begin{proposition}\cite{geom}\label{prop2}
The Levi-Civita connection $\tilde{\nabla}$ of the generalised Cheeger-Gromoll metric $h_{p,q}$ is given by:
\begin{align*}
\tilde{\nabla}_{X^h} Y^h &= (\nabla_{X}Y)^h - \tfrac{1}{2} (R(X,Y)e)^v \\
\tilde{\nabla}_{X^h} Y^v &= (\nabla_{X}Y)^v + \tfrac{\omega^{p}(e)}{2} (R(e,Y)X)^h \\
\tilde{\nabla}_{X^v} Y^h &= \tfrac{\omega^{p}(e)}{2} (R(e,X)Y)^h \\
\tilde{\nabla}_{X^v} Y^v &= -p\omega(e) \left[ g(X,e) Y^v + g(Y,e)X^v \right] +\\ 
& \tfrac{p\omega(e)+q}{1+ q|e|^2}g(X,Y) U + \tfrac{pq\omega(e)}{1+ q|e|^2}g(X,e)g(Y,e) U ,
\end{align*}
where $R$ is the curvature tensor of $(M,g)$ and $U$ the canonical vertical vector field.
\newline Note that $\tilde{\nabla}_{X^v} Y^v \in V$ so the fibres are totally geodesic.
\end{proposition}

This enables us to compute the tension field:

\begin{proposition}
A vector field $\sigma : (M,g) \to (TM,h_{p,q})$ is a harmonic map if and only if:
\begin{equation}
\left\lbrace
\begin{array}{l}
 \sum_{i=1}^n  R(\sigma, \nabla_{e_{i}}\sigma)e_{i} =0 \\
  \nabla^{*}\nabla \sigma + 2p \omega(\sigma) \nabla_{X(\sigma)}\sigma = \bigg[\tfrac{p\omega(\sigma) + q}{1+ 
q|\sigma|^2}|\nabla\sigma|^2 + 
\tfrac{pq\omega(\sigma)}{1+ q|\sigma|^2} |X(\sigma)|^2 \bigg]\sigma
\end{array}
\right.
\end{equation}
where $\{e_{i}\}_{i=1,\dots,n}$ is an orthonormal frame on $(M,g)$.
\end{proposition}

\begin{proof}
A vector field being an immersive bijection, $\sigma(M)\subset TM$ 
is a submanifold of $TM$.\\
Therefore, using a orthonormal frame $\{e_{i}\}_{i=1,\dots,n}$ on $M$:
\begin{align*}
&\tau(\sigma) = \trace \nabla d\sigma =  \sum_{i=1}^{n} \tilde{\nabla}_{d \sigma(e_{i})}d \sigma(e_{i}) - d \sigma (\nabla_{e_{i}} e_{i}) 
\\
&= \sum_{i=1}^{n} \tilde{\nabla}_{(e_{i})^h + (\nabla_{e_{i}}\sigma)^v} [(e_{i})^h + (\nabla_{e_{i}}\sigma)^v] - (\nabla_{e_{i}} e_{i})^h 
- 
(\nabla_{\nabla_{e_{i}} e_{i}}\sigma)^v \\
&= \bigg[ \nabla^{*}\nabla \sigma  -2p \omega(\sigma) 
\nabla_{X(\sigma)}\sigma + \tfrac{p\omega(\sigma) + q}{1+ q|\sigma|^2} |\nabla\sigma|^2\sigma + \tfrac{pq\omega(\sigma)}{1+ q|\sigma|^2} 
|X(\sigma)|^2 \sigma \bigg]^v \\
&+ \sum_{i=1}^{n} \omega^m(e) (R(\sigma, \nabla_{e_{i}}\sigma)e_{i})^h 
\end{align*}
\end{proof}

Since the vector tangent to a variation through sections must be vertical, the Euler-Lagrange equation of the functional $E^{v}$ obtained 
by variations through sections (Equation~\eqref{eq1}) must be equivalent to the vertical part of the tension field. The next 
proposition checks this directly.

\begin{proposition}
A section $\sigma$ of $TM$ satisfies:
\begin{equation}\label{eq6}
(1+|\sigma|^2) \nabla^{*}\nabla \sigma + 2p \nabla_{X(\sigma)}\sigma = \left[p|\nabla\sigma|^{2} -pq |X(\sigma)|^2 
-q(1+|\sigma|^2)\Delta\tfrac{|\sigma|^2}{2} \right] \sigma
\end{equation}
if and only if
\begin{equation}\label{eq7}
\nabla^{*}\nabla \sigma + 2p \omega(\sigma) \nabla_{X(\sigma)}\sigma = \bigg[\tfrac{p\omega(\sigma) + q}{1+ q|\sigma|^2}|\nabla\sigma|^2 
+ \tfrac{pq\omega(\sigma)}{1+ q|\sigma|^2} |X(\sigma)|^2 \bigg]\sigma
\end{equation}
\end{proposition}
\begin{proof}
Assume \eqref{eq6} and $q\neq 0$, then \eqref{eq7} is equivalent to 
\begin{align*}
&  -q|\sigma|^2(1+|\sigma|^2) \nabla^{*}\nabla \sigma -2pq|\sigma|^2\nabla_{X(\sigma)}\sigma + q(1+|\sigma|^2) 
(\Delta\tfrac{|\sigma|^2}{2}) \sigma  \\
&\qquad  \qquad+\bigg[ q(1+|\sigma|^2)|\nabla\sigma|^2 + 2pq |X(\sigma)|^2 \bigg] \sigma = 0\\
&\Leftrightarrow |\sigma|^2 \left( (1+|\sigma|^2) \nabla^{*}\nabla \sigma + 2p\nabla_{X(\sigma)}\sigma\right) = 
\bigg[ (1+|\sigma|^2)(|\nabla\sigma|^2 + \Delta\tfrac{|\sigma|^2}{2}) + 2p |X(\sigma)|^2 \bigg] \sigma \\
&\Leftrightarrow |\sigma|^2 \left( (1+|\sigma|^2) \nabla^{*}\nabla \sigma + 2p\nabla_{X(\sigma)}\sigma\right) =
g( (1+|\sigma|^2)\nabla^{*}\nabla \sigma + 2p \nabla_{X(\sigma)}\sigma, \sigma) \sigma\\
&\mbox{ since } g(\nabla_{X(\sigma)}\sigma, \sigma) = |X(\sigma)|^2 
\end{align*}
Therefore \eqref{eq7} is equivalent to $|\sigma|^2 \Lambda = g(\Lambda,\sigma) \sigma$ where $\Lambda= (1+|\sigma|^2)\nabla^{*}\nabla 
\sigma + 2p \nabla_{X(\sigma)}\sigma$, and this is true if and only if $\Lambda=f\sigma$, which is precisely what \eqref{eq6} says.

Conversely, assume that \eqref{eq7} holds, then \eqref{eq6} is equivalent to 
\begin{align*}
& (1+|\sigma|^2) \nabla^{*}\nabla \sigma + 2p \nabla_{X(\sigma)}\sigma = \\
& \qquad \qquad \bigg[ p|\nabla \sigma|^2 - pq |X(\sigma)|^2 - q(1+|\sigma|^2) ( g(\nabla^{*}\nabla \sigma,\sigma) - |\nabla \sigma|^2) 
\bigg] \sigma \\
&\Leftrightarrow (1+|\sigma|^2) \nabla^{*}\nabla \sigma + 2p \nabla_{X(\sigma)}\sigma + q \left[ g((1+|\sigma|^2) \nabla^{*}\nabla \sigma 
+ p \nabla_{X(\sigma)}\sigma,\sigma)\right] \sigma  \\
& \qquad \qquad = (p+q +q|\sigma|^2)|\nabla \sigma|^2 \sigma \\
&\Leftrightarrow \left[ (1+|\sigma|^2) \nabla^{*}\nabla \sigma + 2p \nabla_{X(\sigma)}\sigma\right] (1+q|\sigma|^2) =\\
& \qquad \qquad \left[(p+q +q|\sigma|^2)|\nabla \sigma|^2 + pq|X(\sigma)|^2\right] \sigma
\end{align*}
since, again, \eqref{eq7} implies that $|\sigma|^2 \Lambda = g(\Lambda,\sigma) \sigma$.
\end{proof}

\begin{proposition}\label{prop5}
The scaled Hopf vector fields are harmonic maps from ${\sn}^{2n+1}$ to $T{\sn}^{2n+1}$ equipped 
with the generalised Cheeger-Gromoll metric $h_{p,q}$ ($p>1$, $q\geqslant0$).
\end{proposition}
\begin{proof}
On spheres, the horizontal part of the tension field becomes:
\begin{align*}
\sum_{i=1}^{n} R(\sigma, \nabla_{e_{i}} \sigma ) e_{i} &= (\Div{\sigma}) \sigma - \nabla_{\sigma}\sigma
\end{align*}
So if $\sigma$ is Killing of constant norm, it is harmonic.
\end{proof}

\begin{remark}
a) For $2\geqslant p >1$ and $q=0$, these represent new examples of harmonic maps 
from a compact manifold into positive sectional curvature (cf.~\cite{geom}).\\
b) No Killing vector field on ${\sn}^{2n}$ can be a harmonic map.
\end{remark}

\begin{example}
Proposition~\ref{prop5} also extends to the vector fields on Heisenberg and $\widetilde{SL_{2}(\rn)}$ 
considered in Example~\ref{example1}:\\
i) $H_{3}$: With the same notations, let $\sigma$ be the vector $E_{3}$, since the non-zero components of the curvature tensor of $H_3$
are $R_{1221}=\tfrac{3}{4} \mbox{ and } R_{1331}=-\tfrac{1}{4}=R_{2332},$
then $\sum_{i=1}^{3} R(\sigma, \nabla_{e_{i}} \sigma ) e_{i}  = 0$.\\
ii) $\widetilde{SL_{2}(\rn)}$: Still with the same notations as in Example~\ref{example1}, let $\sigma$ be the vector $E_{1}$, then:
$\sum_{i=1}^{3} R(\sigma, \nabla_{e_{i}} \sigma ) e_{i} = 0$, since $R(E_{1}, E_{3} ) E_{2} = R(E_{1}, E_{2} ) E_{3} = 0$.
\end{example}


\begin{thebibliography}{999}

\bibitem{LBW}
M.~Benyounes, E.~Loubeau and C.~M.~Wood,
\newblock Harmonic sections of Riemannian vector bundles and metrics of Cheeger-Gromoll type,
\newblock to appear in {\em Diff. Geom. Appl.}

\bibitem{geom}
M.~Benyounes, E.~Loubeau and C.~M.~Wood,
\newblock The geometry of generalised Cheeger-Gromoll metrics,
\newblock in preparation.

\bibitem{E-L}
J.~Eells and L.~Lemaire,
\newblock {\em Selected topics on harmonic maps},
\newblock AMS, 1983.

\bibitem{Gud-Kap-Expo}
S.~Gudmundsson and E.~Kappos,
\newblock On the geometry of tangent bundles,
\newblock {\em Expositiones Mathematicae}, 20 (2002), 1--14.

\bibitem{Gud-Kap-Tok}
S.~Gudmundsson and E.~Kappos,
\newblock On the geometry of tangent bundle with the Cheeger-Gromoll metric,
\newblock {\em Tokyo J. Math.}, 25 (2002), 75--83.

\bibitem{Kowal}
O.~Kowalski
\newblock Curvature of the induced Riemannian metric on the tangent bundle of a Riemannian manifold,
\newblock {\em J. Reine Angew. Math.}, 250 (1971) 124--129. 

\bibitem{Cez}
C.~Oniciuc,
\newblock The tangent bundle and harmonicity,
\newblock {\em An. Stiint. Univ. Al. I. Cuza Mat.}, 43(1997), 151--172.

\bibitem{Petersen}
P.~Petersen,
\newblock {\em Riemannian geometry},
\newblock Springer-Verlag, 1998.

\bibitem{Yau}
S.~T.~Yau,
\newblock Some function-theoretic properties of complete Riemannian manifolds and their applications to geometry,
\newblock {\em Indiana M. J.} 25 (1976), 659--670.

\end{thebibliography}
\end{document}